\documentclass[11pt]{amsart}
\usepackage{amssymb}
\usepackage{dsfont}
\usepackage{color}

\newtheorem{theorem}{Theorem}
\newtheorem{lemma}[theorem]{Lemma}
\newtheorem{proposition}[theorem]{Proposition}
\newtheorem{problem}[theorem]{Problem}
\newtheorem{corollary}[theorem]{Corollary}

\theoremstyle{definition}
\newtheorem{df}{Definition}

\newtheorem{remark}[df]{Remark}
\newcommand{\B}{\mathcal B}

\newcommand{\N}{\mathbb N}

\newcommand{\R}{\mathbb R}

\newcommand{\on}{\operatorname}
\newcommand{\Si}{\Sigma}
\newcommand{\mc}{\mathcal}

\author{Marek Balcerzak}
\address{Institute of Mathematics, Lodz University of Technology, al. Politechniki 8, 93-590
\L\'od\'z, Poland}
\email{marek.balcerzak@p.lodz.pl}
\author{Jacek Hejduk}
\address{Faculty of Mathematics and Computer Science, University of Lodz, ul. Banacha 22,
90-238 \L\'od\'z, Poland}
\email {jacek.hejduk@wmii.uni.lodz.pl}
\author{Artur Wachowicz}
\address{Institute of Mathematics, Lodz University of Technology, al. Politechniki 8, 93-590
\L\'od\'z, Poland}
\email{artur.wachowicz@p.lodz.pl}

\title{Baire category lower density operators with Borel values}
\subjclass[2010]{Primary: 28A51; Secondary: 03E50, 03E35}
\keywords{lower density operator, density point, the Baire property, meager set, the Cantor space}
\date{}

\begin{document}

\begin{abstract}
We prove that the lower density operator associated with the Baire category density points in the real line has Borel values of class $\pmb \Pi^0_3$ which is analogous to the measure case. We also introduce the notion of the Baire category density point of a subset with the Baire property of the Cantor space, and we prove that it generates a lower density operator with Borel values of class $\pmb \Pi^0_3$.
\end{abstract}
\maketitle

\section{Introduction}
Let $\Si$ be a $\sigma$-algebra of subsets of a set $X\neq\emptyset$ and let $J\subset \Si$ be a 
$\sigma$-ideal.
We write $A\sim B$ whenever the symmetric difference $A\bigtriangleup B$ belongs to $J$.
Note that $\sim$ is the equivalence relation on $\Si$.
A mapping
$\Phi\colon \Si\to \Si$ is called a {\em lower density operator} with respect to $J$ if the following conditions are satisfied:
\begin{itemize}
\item[(i)] $\Phi(X)=X$ and $\Phi(\emptyset)=\emptyset$,
\item[(ii)] $A\sim B\implies\Phi(A)=\Phi(B)$ for every $A,B\in \Si$,
\item[(iii)] $A\sim\Phi(A)$ for every $A\in \Si$,
\item[(iv)] $\Phi(A\cap B)=\Phi(A)\cap\Phi(B)$ for every $A,B\in \Si$.
\end{itemize}
The book \cite{Ox} gives a standard example of a lower density operator for the $\sigma$-algebra $\mathcal L$ of Lebesgue measurable subsets of $\R$, with respect to the $\sigma$-ideal $\mathcal N$ of null sets. This operator $\Phi\colon\mathcal L\to\mathcal L$ assigns to $A\in\mathcal L$ its set of density points, that is
$$\Phi(A):=\left\{ x\in\R\colon\lim_{h\to 0^+}\frac{\lambda(A\cap [x-h,x+h])}{2h}=1\right\} $$
where $\lambda$ denotes Lebesgue measure. It is known (cf. \cite{W}) that the value $\Phi(A)$ is always 
a Borel set of type $F_{\sigma\delta}$ (or $\pmb \Pi^0_3$ in the modern notation; cf. \cite{Ke}).
There is an exact counterpart of this operator in the case when $\R$ is replaced by 
the Cantor space $\{0,1\}^\N$ with the respective product measure (see \cite{AC}; 
cf. also \cite[Exercise 17.9]{Ke}). It was proved in \cite{AC} that the values $\Phi(A)$ in this case are again in the Borel class $\pmb \Pi^0_3$, and this Borel level cannot be lower for a large class of measurable sets. For further results in this case, see \cite{AC1}, \cite{ACC}.

The Baire category analogue of the notion of a density point is due to Wilczy\'nski \cite{W1}. Several further results are provided in \cite{PWW}. This notion led to the lower density operator $\Phi_\mathcal M$ for the $\sigma$-algebra $\mathcal B$ of subsets of $\R$ having the Baire property, with respect to the $\sigma$-ideal $\mathcal M$ of meager subsets of $\R$.

Let us recall the respective definitions. Let $A\in\mathcal B$ and $x\in\R$. We say that:
\begin{itemize}
\item $0$ is an $\mathcal M$-density point of $A$ if, for each increasing sequence $(n_k)$ of positive integers, there is a subsequence $(n_{i_k})$ such that 
$$\limsup_{k\in\N}((-1,1)\setminus n_{i_k} A)\in\mathcal M$$
where $\alpha A:=\{ \alpha t\colon t\in A\}$ for $\alpha\in\R$;
\item $x$ is an $\mathcal M$-density point of $A$ if $0$ is an $\mathcal M$-density point of $A-x:=
\{t-x\colon t\in A\}$;
\item $x$ is an $\mathcal M$-dispersion point of $A$ if it is an $\mathcal M$-density point of $A^c:=\R\setminus A$. So, $0$ is an $\mathcal M$-dispersion point of $A$ whenever, for each increasing sequence $(n_k)$ of positive integers, there is a subsequence $(n_{i_k})$ such that 
$$\limsup_{k\in\N}((-1,1)\cap n_{i_k} A)\in\mathcal M.$$
\end{itemize}
Note that, if we replace $\mathcal B$ by $\mathcal L$, and $\mathcal M$ by $\mathcal N$ in the obove definitions,
we obtain the classical notions of density and dispersion points; see \cite{PWW} and \cite[p. 7]{CLO}.
The operator $\Phi_\mathcal M\colon \mathcal B\to\mathcal B$ that assigns to $A\in\mathcal B$ its set of $\mathcal M$-density points, satisfies conditions (i)--(iv); cf. \cite{PWW} and \cite[Lemma 2.2.1]{CLO}.

Our purposes in this paper are twofold. Firstly, we prove (in Section~2) that the values of the operator 
$\Phi_{\mc M}$ are Borel of bounded level. Namely, they hit into class $\pmb \Pi^0_3$ which is analogous to the measure case. We use, as the main main tool, a combinatorial characterization of 
$\Phi_{\mc M}(G)$, for an open set $G$, based on ideas of E. {\L}azarow \cite{EL}. This characterization has motivated us to define 
(in Section~3) the notion of an $\mc M$-density point of a set with the Baire property in the Cantor space. Then we show that the respective mapping $\Phi_{\mc M}$ is a lower density operator, which fulfills our second purpose. The final Section~4 contains additional remarks and open problems.

\section{The values of $\Phi_M$ are Borel}
It is well known (cf. \cite[Thm 4.6]{Ox}) that a set $A\subseteq\R$ with the Baire property can be expressed in the form $A=G\bigtriangleup E$ where $G$ is open and $E\in\mathcal M$. Moreover, such an expression is unique, if we additionally assume that $G$ is regular open, that is 
$G=\on{int}(\on{cl} G)$.
This regular open set will be denoted by $\widetilde{A}$ and called the regular open kernel of $A$.
Then $\widetilde{A}$ is the largest, in the sense of inclusion, among open sets in the above expression.
We have $A\sim\widetilde{A}$, and so, $\Phi_\mathcal M(A)=\Phi_\mathcal M(\widetilde{A})$. 

We will use the following characterization which is due to E. {\L}azarow \cite{EL}.
We reformulate a bit the original version and give the proof for the reader's convenience. However,
we strongly mimic the idea from \cite{EL}. For a related characterization, see \cite[Thm 2.2.2] {CLO}.

\begin{proposition}\cite{EL} \label{pp}
The number $0$ is an $\mathcal M$-dispersion point of an open set $G\subseteq\R$ if and only if the following holds:
\begin{quotation}
for each $n\in\N$ there exists $k\in\N$ such that, for all integer $\ell >k$ and every $i\in\{ -n,\dots, n-1\}$, there exists $j\in\{1,\dots ,k\}$ satisfying
$$G\cap\left(\frac{1}{\ell}\left(\frac{i}{n}+\frac{j-1}{nk}\right),\;\frac{1}{\ell}\left(\frac{i}{n}+\frac{j}{nk}\right)\right)=\emptyset .$$
\end{quotation}
\end{proposition}
\begin{proof} Necessity. Suppose that it is not the case. Then there exists $n_0\in\N$ such that, for each $k\in\N$, we can pick integers $\ell_k>k$ and $i_k\in\{-n_0,\dots , n_0-1\}$ such that, for each $j\in\{ 1,\dots , k\}$,
$$G\cap\left(\frac{1}{\ell_k}\left(\frac{i_k}{n_0}+\frac{j-1}{n_0k}\right),\;\frac{1}{\ell_k}\left(\frac{i_k}{n_0}+\frac{j}{n_0k}\right)\right)\neq\emptyset .$$
We may assume that $\ell_{k+1}>\ell_k$ for every $k$. Pick a value of the sequence $(i_k)$
in $\{-n_0,\dots ,n_0-1\}$ that appears in the above condition infinitely many times.
We call it $i_0$, and let it be associated with a subsequence $(\ell_{k_m})$ of $(\ell_k)$.
Take an arbitrary subsequence $(r_m)$ of $(\ell_{k_m})$. It follows that, for every $p\in\N$, the set
$$\bigcup_{m\ge p}(r_m G)\cap\left(\frac{i_0}{n_0},\frac{i_0+1}{n_0}\right)$$
is open and dense in $\left(\frac{i_0}{n_0},\frac{i_0+1}{n_0}\right)$. Hence the set
$$\bigcap_{p\in\N}\bigcup_{m\ge p}(r_m G)\cap\left(\frac{i_0}{n_0},\frac{i_0+1}{n_0}\right)$$
is comeager in $\left(\frac{i_0}{n_0},\frac{i_0+1}{n_0}\right)$. Consequently,
$$\limsup_{m\in\N}\left((-1,1)\cap r_m G\right)\supseteq \limsup_{m\in\N}\left(\left(\frac{i_0}{n_0},\frac{i_0+1}{n_0}\right)\cap( r_m G)\right)\notin\mathcal M.$$
It shows that $0$ is not an $\mathcal M$-dispersion point of $G$. Contradiction.

Sufficiency. Let $(m_s)$ be an increasing sequence of positive integers. We will define inductively
the respective subsequence $(r_s)$ of $(m_s)$ which witnesses that $0$ is an $\mathcal M$-dispersion point of $A$. In fact, we will define the family $\{S^{(n)}\colon n\in\N\}$ of subsequences of $(m_s)$ where $S^{(n+1)}$ is a subsequence of $S^{(n)}$ for every $n$. Taking the first term from every $S^{(n)}$, we obtain the sequence $(r_s)$.

By the assumption, for $n:=1$ pick $k_1\in\N$ such that for all integers $\ell> k_1$
and $i\in\{-1,0\}$ we can choose $j=j(\ell,i)\in\{1,\dots , k_1\}$ with
$$G\cap\left(\frac{1}{\ell}\left({i}+\frac{j-1}{k_1}\right),\;\frac{1}{\ell}\left({i}+\frac{j}{k_1}\right)\right)=\emptyset .$$
Firstly, for $i:=-1$ we find a subsequence $S$ of $(m_s)$ such that we obtain the same value $j_{-1}$ of 
$j(\ell, -1)$ for all $\ell>k_1$ from $S$, and then we find a subsequence $S^{(1)}$ of $S$
such that we obtain the same value $j_{0}$ of $j(\ell,0)$ for all $\ell$ from $S^{(1)}$.
In this way, we obtain a subsequence $(m_{\alpha_1(s)})=:S^{(1)}$ in the first step of induction.
Let $r_1:=m_{\alpha_1(1)}$. We may assume that $r_1>k_1$.

Now, let $n>1$. Suppose that we have $k_{n-1}\in\N$ and  
a subsequence $(m_{\alpha_{n-1}(s)})=:S^{(n-1)}$
of $(m_s)$ such that $r_{n-1}:=m_{\alpha_{n-1}(1)}>k_{n-1}$ and the following is true
\begin{quotation}
 for each $i\in\{-n+1,\dots ,n-2\}$
an integer $j_i\in\{1,\dots ,k_{n-1}\}$ is chosen so that for all $\ell$ in $S^{(n-1)}$,
$$
G\cap\left(\frac{1}{\ell}\left(\frac{i}{n-1}+\frac{j_i-1}{(n-1)k_{n-1}}\right),\;\frac{1}{\ell}\left(\frac{i}{n-1}+\frac{j_i}{(n-1)k_{n-1}}\right)\right)=\emptyset .
$$
\end{quotation}
By the assumption, for the given $n$, pick $k_n\in\N$ such that, for all $\ell> k_n$
and $i\in\{-n,\dots n-1\}$, there exists $j=j(\ell,i)\in\{1,\dots , k_n\}$ with 
$$G\cap\left(\frac{1}{\ell}\left(\frac{i}{n}+\frac{j-1}{nk_n}\right),\;\frac{1}{\ell}\left(\frac{i}{n}+\frac{j}{nk_n}\right)\right)=\emptyset .$$
We may assume that $k_n>k_{n-1}$. Note that, if $\ell>k_n$ is in the subseqence $S^{(n-1)}$, then
for every $i\in\{-n+1,\dots ,n-2\}$, the choice $j=j_i$ (independent of $\ell$) made in the previous step will be preserved. Next, as in the first step, after a double choice, we pick a subsequence $S^{(n)}=(m_{\alpha_{n}(s)})$
of $S^{(n-1)}$ 
such that, for each $i\in \{-n, n-1\}$, the value $j_i$ of $j(\ell,i)$ is the same whenever $\ell>k_n$ is taken from $S^{(n)}$. Also, we may assume that $m_{\alpha_{n}(1)}>\max\{ k_n, m_{\alpha_{n-1}(1)}\}$. Finally, let $r_{n}:=m_{\alpha_{n}(1)}$. This ends the construction.

The above construction guarantees that, if $n_0\in\N$ is fixed and $n>n_0$, then for all terms $\ell$ in $S^{(n)}$ and every
$i\in\{-n_0,\dots, n_0-1\}$, the choice $j=j_i$ made in the $n_0$-th step of construction remains
unchanged in the $n$-th step.

Clearly, $(r_n)$ defined in the construction is a subsequence of $(m_n)$. We will show that the set
$$\limsup_{n\in\N}((-1,1)\cap r_n G)$$
 is nowhere dense which implies that $0$ is an $\mathcal M$-dispersion point of $G$. Let $(a,b)$ be any subinterval of $(-1,1)$. Pick $n_0\in\N$ and $i_0\in\{-n_0,\dots ,n_0-1\}$ such that $\left(\frac{i_0}{n_0}, \frac{i_0+1}{n_0}\right)\subseteq (a,b)$.
By the construction and the above remark, for each $n\ge n_0$, the integer $r_{n}$ is in the sequence $(m_{\alpha_{n_0}(s)})$,
and there exists $j\in\{1,\dots ,k_{n_0}\}$ such that for each $n\ge n_0$,
$$G\cap\left(\frac{1}{r_{n}}\left(\frac{i_0}{n_0}+\frac{j-1}{n_0k_{n_0}}\right),\;\frac{1}{r_{n}}\left(\frac{i_0}{n_0}+\frac{j}{n_0k_{n_0}}\right)\right)=\emptyset .$$
Put $(a_0,b_0):=\left(\frac{i_0}{n_0}+\frac{j-1}{n_0k_{n_0}},\;\frac{i_0}{n_0}+\frac{j}{n_0k_{n_0}}\right).$ Then $(a_0,b_0)\subseteq (a,b)$ and
$(a_0,b_0)\cap\bigcup_{n\ge n_0}(r_n G)=\emptyset$.
Consequently, $$(a_0,b_0)\subseteq(a,b)\setminus \limsup_{n\in\N}((-1,1)\cap r_n G)$$
which means that $\limsup_{n\in\N}((-1,1)\cap r_n G)$ is nowhere dense,
as desired.
\end{proof}

Proposition \ref{pp} can be easily reformulated in the case where $0$ is replaced by a point $x\in\R$.
Furthermore, $x$ is an $\mc M$-density point of $A\in\mc B$ if and only if
$0$ is an $\mc M$-dispersion point of $G-x$ where
$A^c=G\bigtriangleup E$, $G$ is open and $E\in\mc M$. In particular, we can take the regular open kernel of $A^c$ in the role of $G$. So, we obtain the following

\begin{corollary} \label{cc}
A point $x\in\R$ is an $\mathcal M$-density point of a set $A\in\B$ if and only if
$$(\forall\; n\in\N)(\exists\; k\in\N)(\forall\; \ell>k)(\forall\; i\in\{ -n,\dots ,n-1\})(\exists\; j\in\{ 1,\dots ,k\})$$
$$(\widetilde{A^c}-x)\cap\left(\frac{1}{\ell}\left(\frac{i}{n}+\frac{j-1}{nk}\right),\;\frac{1}{\ell}\left(\frac{i}{n}+\frac{j}{nk}\right)\right)=\emptyset .$$
\end{corollary}

\begin{lemma} \label{LL}
Given an open set $U\subseteq \R$ and an interval $(a,b)\subseteq (-1,1)$, the set
$$E:=\{ x\in\R\colon (U-x)\cap (a,b)=\emptyset\}$$ is closed.
\end{lemma}
\begin{proof}
Note that $E^c=\pi_1[D]$ where
$$D:=\{(x,t)\in\R\times (-1,1)\colon t+x\in U\;\mbox{ and }t\in(a,b)\}$$ and $\pi_1\colon\R\times [-1,1]\to\R$ is the projection $\pi_1(x,t):=x$. Clearly, the set $D$ is open. 
Hence $\pi_1[D]$ is open and consequently, $E$ is closed.
\end{proof}

\begin{theorem} \label{TT}
For each $A\in\mathcal B$, the set $\Phi_{\mathcal M}(A)$, of all $\mathcal M$-density points of $A$,
is a Borel set of type $F_{\sigma\delta}$, i.e. of class $\pmb \Pi_3^0$.
\end{theorem}
\begin{proof}
Let $A\in\mathcal B$. Fix $n,k\in\N$, $\ell>k$, $i\in\{-n,\dots ,n-1\}$ and 
$j\in\{1,\dots ,k\}$. Denote
$$V_{n,k,\ell,i,j}:=\left(\frac{1}{\ell}\left(\frac{i}{n}+\frac{j-1}{nk}\right),\;\frac{1}{\ell}\left(\frac{i}{n}+\frac{j}{nk}\right)\right).$$
Put $E:=\{ x\in\R\colon(\widetilde{A^c}-x)\cap V_{n,k,\ell,i,j}=\emptyset\}$.
By Lemma \ref{LL} this set is closed.
This together with Corollary~\ref{cc} gives the assertion.
\end{proof}

\begin{remark} In the papers \cite{W1}, \cite{PWW}, the ideal of meager sets in $\R$ is denoted by $I$.
The respective topology defined by the operator $\Phi_I$($=\Phi_{\mc M}$) is the Baire category analogue of the density topology in $\R$ and is called the $I$-density topology. 
Characterizations similar to that given in Proposition \ref{pp} were applied to the so-called 
$I$-approximate derivatives. In \cite{EL} it was proved that $I$-approximate 
derivative is of Baire class 1. Another theorem on $I$-approximate derivatives was obtained in \cite{BW}
where the characterization from \cite[thm 2.2.2 (vii)]{CLO} was used. In fact, the characterization from \cite{EL} turns out more suitable in that case which was shown in the PhD thesis \cite{Wa} of the third author.
\end{remark}

\section{$\mathcal M$-density points in the Cantor space}
We will use the characterization stated in Proposition \ref{pp} to introduce the notion
of an $\mathcal M$-density point of a set with the Baire property in the Cantor space $\{0,1\}^\N$. 
The families of meager sets and of sets with the Baire property in $\{0,1\}^\N$ will be denoted again by $\mathcal M$ and $\mathcal B$, respectively. 
Again, every set $A\in\mathcal B$ can be uniquelly expressed in the form $A=G\bigtriangleup E$ where $G$ is regular open and $E\in\mathcal M$ (cf. \cite [Exercise 8.30]{Ke}). Then $G$ will be denoted by $\widetilde A$ and called the regular open kernel of $A$.

Recall that sets of the base in the product topology of the Cantor space are of the form
$$U(s):=\{ x\in\{0,1\}^\N\colon s\subseteq x\}$$
for any finite sequence $s$ of zeros and ones (that is, $s\in\{0,1\}^{<\N}$).
Given $x\in \{0,1\}^\N$, by $x|n$ we denote the restriction of $x$ to the first $n$ terms.
For $s,t\in\{0,1\}^{<\N}$, let $s\frown t$ denote their concatenation where terms of $t$ follow the terms of $s$.

We say that $x\in\{0,1\}^\N$ is an $\mathcal M$-dispersion point of a set $A\in\mathcal B$
if 
\begin{quotation}
for each $n\in\N$ there exists $k\in\N$ such that, for each $\ell\in\N$ with $\ell>k$, and every $s\in\{0,1\}^n$, there exists $t\in\{0,1\}^k$ with
$$\widetilde{A}\cap U((x|\ell)\frown s\frown t)=\emptyset.$$
\end{quotation}
We say that $x\in\{0,1\}^\N$ is an $\mathcal M$-density point of a set $A\in\mathcal B$
if it is an $\mathcal M$-dispersion point of $A^c$. So, in the above condition, one should replace
$\widetilde{A}$ by $\widetilde{A^c}$. 

Note that we can use the group structure of $\{0,1\}^\N$
(with coordinatewise addition mod~2), and thus
condition $\widetilde{A}\cap U((x|\ell)\frown s\frown t)=\emptyset$ can be written as
$(\widetilde{A}-x)\cap U(({\mathbf 0}|\ell)\frown s\frown t)=\emptyset$ where
$\mathbf 0:=(0,0,\dots )$. We will use this fact in the proof of Theorem \ref{Can}.

\begin{theorem} \label{Can}
Let $\Phi_\mathcal M$ assign to each set $A\subseteq\{ 0,1\}^\N$ with the Baire property, the set of 
$\mathcal M$-density points of $A$. Then $\Phi_\mathcal M\colon\mathcal B\to\mathcal B$ is a lower density operator with respect to $\mathcal M$ and its values are $F_{\sigma\delta}$ sets, i.e. sets
of class $\pmb\Pi^0_3$.
\end{theorem}

\begin{proof}
The final assertion follows from the definition of $\Phi_\mc M$ and the arguments anologous to those used for Theorem \ref{TT}. Let us sketch this proof. For fixed $n,k\in\N$, $\ell>k$, and 
$s\in\{0,1\}^n$, $t\in\{0,1\}^k$, we denote
$V_{n,k,\ell,s,t}:= U(({\mathbf 0}|\ell)\frown s\frown t)$. Then, as in Lemma \ref{LL}, we observe
that the set $\{ x\in\{0,1\}^\N\colon(\widetilde{A^c}-x)\cap V_{n,k,\ell,s,t}=\emptyset\}$
is closed. Finally, it suffices to use the respectively modified condition stating that 
$x\in\Phi_{\mc M}(A)$.

Now, it is enough to check conditions (i)--(iv) stated in the definition of a lower density operator.

Condition (i) is clearly valid. To show (ii) note that the following is true for any $A,B\in\mc B$:
$$A\sim B\implies A^c\sim B^c\implies \widetilde{A^c}= \widetilde{B^c}.$$
This yields (ii) by the definition of an $\mc M$-density point.

Let us prove an additional property. Let $A,B\in\mc B$ and $A\subseteq B$.
We show that $\Phi_\mc M(A)\subseteq \Phi_\mc M(B)$.
From $A\subseteq B$ it follows that $B^c\subseteq A^c$ and then 
$\widetilde{B^c}\subseteq\widetilde{A^c}$.  
Hence $\Phi_\mc M(A)\subseteq \Phi_\mc M(B)$ by the definition of an $\mc M$-density point.

Next, we will prove (iv). Let $A,B\in\mc B$. By the above property we obtain
$$\Phi_\mc M(A\cap B)\subseteq \Phi_\mc M(A)\cap\Phi_\mc M(B).$$
To show the reverse inclusion, let $x\in\Phi_\mc M(A)\cap\Phi_\mc M(B)$.
Fix $n\in\N$. Since $x$ is an $\mc M$-density point of $A$, we pick $k'\in\N$
such that for all $\ell>k'$ and $s\in\{0,1\}^n$, we can find $t'\in\{0,1\}^{k'}$ with
$$\widetilde{A^c}\cap U((x|\ell)\frown s\frown t')=\emptyset.$$
Since $x$ is an $\mc M$-density point of $B$, we consider $n+k'$ instead of $n$ and pick $k''\in\N$ 
such that for all $\ell>k''$ and $s'\in\{0,1\}^{n+k'}$, we can find $t''\in\{0,1\}^{k''}$ with
$$\widetilde{B^c}\cap U((x|\ell)\frown s'\frown t'')=\emptyset.$$
Now, fix $\ell>k'+k''$ and $s\in\{0,1\}^n$. Then taking $s':=s\frown t'$, we pick
$t''\in\{0,1\}^{k''}$ with
$\widetilde{B^c}\cap U((x|\ell)\frown s'\frown t'')=\emptyset$.
Observe that
\begin{equation} \label{raz}
(\widetilde{A^c}\cup\widetilde{B^c})\cap U((x|\ell)\frown s\frown t'\frown t'')=\emptyset.
\end{equation}
Then $k:=k'+k''$ and $t:=t'\frown t''\in\{0,1\}^{k}$ will witness that
$x\in\Phi_\mc M(A\cap B)$ provided that $\widetilde{A^c}\cup\widetilde{B^c}$ can be replaced by
$\widetilde{D}$ for $D:=(A\cap B)^c$ in condition (\ref{raz}).
So, let us show this last requirement. First note that $\widetilde{A^c}\cup\widetilde{B^c}$
can be replaced by $C:=\on{int}\on{cl}(\widetilde{A^c}\cup\widetilde{B^c})$ in condition (\ref{raz})
(since the set $U(\cdot )$ is open).
Additionally, $C\sim \widetilde{A^c}\cup\widetilde{B^c}$. Also, from $\widetilde{A^c}\sim A^c$ and
$\widetilde{B^c}\sim B^c$ it follows that 
$C\sim \widetilde{A^c}\cup\widetilde{B^c}\sim A^c\cup B^c= (A\cap B)^c=D\sim\widetilde{D}$. Since $C$ is a regular open
set (cf. \cite[p. 23]{JW}), we have
$C=\widetilde{D}$ which yields the required condition.

To prove (iii) we need the following property
\begin{equation} \label{dwa}
G\subseteq\Phi_{\mc M}(G)\subseteq\on{cl}(G)\quad
\mbox{ for every open set } G\subseteq\{0,1\}^\N .
\end{equation}
So let $G\subseteq\{0,1\}^\N$ be nonempty open. Let $x\in G$. Fix $n,k\in\N$ and pick $\ell >k$ such that $U(x|\ell)\subseteq G$. Then for all $s\in\{0,1\}^n$ and
$t\in\{0,1\}^k$ we have $U:=U((x|\ell)\frown s\frown t)\subseteq G$. Hence $U\cap G^c=\emptyset$. Since $G^c$ is closed, we have $\widetilde{G^c}\subseteq G^c$ and so, $U\cap \widetilde{G^c}=\emptyset$. 
Thus $x\in\Phi_{\mc M}(G)$. This yields the first inclusion in (\ref{dwa}).
Now, let $x\in\Phi_{\mc M}(G)$ and suppose that $x\notin\on{cl}(G)$. Thus $x$ belongs to the open set $V:=(\on{cl}(G))^c$.
Using the inclusion proved before, we have $x\in\Phi_{\mc M}(V)$. Hence by (i) and (iv),
$$x\in\Phi_{\mc M}(G)\cap\Phi_{\mc M}(V)=\Phi_{\mc M}(G\cap V)=\Phi_{\mc M}(\emptyset)=\emptyset.$$
Contradiction.
This ends the proof of (\ref{dwa}).

Now, fix $A\in\mc B$. Taking $G:=\widetilde{A}$, by (\ref{dwa}) we have 
$\widetilde{A}\sim\Phi_{\mc M}(\widetilde{A})$. Then using (ii) we obtain
$$A\sim\widetilde{A}\sim\Phi_{\mc M}(\widetilde{A})=\Phi_{\mc M}(A)$$
which yields (iii).
\end{proof}

Note that the lower density operator $\Phi_{\mc M}$ generates a topology
$${\mc T}_{\mc M}:=\{ A\in\mc B\colon A\subseteq\Phi_{\mc M}(A)\}$$
that is finer than the standard topology in $\{0,1\}^\N$.
Indeed, every open set in the standard topology belongs to $\mc T_{\mc M}$ by (\ref{dwa}).
Every nontrivial comeager set in the standard topology witnesses that these two topologies are different.
For other properties and their proofs, see \cite{PWW} or \cite[Sec. 2.3]{CLO} where the analogous topology in $\R$ was investigated.

\section{Final remarks}
The presented results should initiate further studies.
We hope that our notion of an $\mc M$-density point is a good Baire category counterpart
of an density point for the respective subsets of the Cantor space.
In particular, an interesting question appears whether the results analogous to those obtained 
in \cite{AC} can be proved.

\begin{problem}
Find natural examples of sets  $A\subseteq \{0,1\}^\N$ (for instance, open or closed)
such that $\Phi_{\mc M}(A)$ is complete $\pmb\Pi^0_3$.
\end{problem}

The Baire category analogue of a density point for subsets with the Baire property of $\R$,
due to Wilczy\'nski, was well motivated by the respective characterization in the measure case where
convergence in measure of characteristic functions and other properties were considered (see \cite{PWW}
or \cite{CLO}). This idea can be used for other $\sigma$-algebras and $\sigma$-ideals in the Euclidean spaces. Such a process was successful in the article \cite{BH} where the lower density operators
associated with the product ideals $\mc M\otimes\mc N$ and $\mc N\otimes\mc M$ in $\R^2$
were defined.

\begin{problem}
Are the lower density operators,
associated with the product ideals $\mc M\otimes\mc N$ and $\mc N\otimes\mc M$ in $\R^2$,
Borel-valued?
\end{problem}
 
 Quite recently, Wilczy\'nski  proposed in \cite{WW} another natural Baire category notion of a density point, called an intensity point, for subsets with the Baire property in $\R$. Then the respective mapping $\Phi_{i}\colon\mc B\to\mc B$ is a lower density operator which produces a topology non-homeomorphic to the $I$-density topology (for $I:=\mc M$) in $\R$.
 
 \begin{problem}
Is the lower density operator $\Phi_i$ Borel-valued? Can one define its analogue in the Cantor space
setting?
\end{problem}


\begin{thebibliography}{a,b,c}
\bibitem{AC} A. Andretta, R. Camerlo, {\it The descriptive set theory of the Lebesgue Density Theorem}, Adv. Math. {\bf 244} (2013), 1--42.
\bibitem{AC1} A. Andretta, R. Camerlo, {\it Analytic sets of reals and the density function in the Cantor space}, European J. Math. {\bf 5} (2019), 49--80.
\bibitem{ACC} A. Andretta, R. Camerlo, C. Costantini, {\it Lebesgue density and exceptional points}, Proc. London Math. Soc. {\bf 118} (2019), 103--142.
\bibitem {BH} M. Balcerzak, J. Hejduk, {\it Density topologies for products of $\sigma$-ideals}, Real Anal. Exchange {\bf 20} (1994-95), 163--177.
\bibitem{BW} M. Balcerzak, A. Wachowicz, {\it Some examples of meager sets in Banach spaces}, Real Anal. Exchange {\bf 26} (2000-2001), 877--884.
\bibitem{CLO} K. Ciesielski, L. Larson, K. Ostaszewski, {\it $\mathcal I$-Density Continuous Functions},
Mem. Amer. Math. Soc. {\bf 107}, no.~515, (1994).
\bibitem{JW} W. Just, M. Weese, {\it Discovering Modern Set Theory II}, Amer. Math. Soc., Providence 1997.
\bibitem{Ke} A.S. Kechris, {\it Classical Descriptive Set Theory}, Springer, New York, 1995.
\bibitem{EL} E. {\L}azarow, {\it On the Baire class of $\mathcal I$-approximate derivatives}, Proc. Amer. Math. Soc. {\bf 100} (1987), 669--674.
\bibitem{Ox} J.C. Oxtoby, {\it Measure and Category}, Springer, New York, 1980.
\bibitem{PWW} W. Poreda, E. Wagner-Bojakowska, W. Wilczy\'nski, {\it A category analogue of the density topology}, Fund. Math. {\bf 125} (1985), 167--173.
\bibitem{Wa} A. Wachowicz, {\it On some residual sets}, PhD thesis, Lodz University of Technology, 
{\L}\'od\'z 2004 (in Polish).
\bibitem{W1} W. Wilczy\'nski, {\it A generalization of the density topology}, Real Anal. Exchange {\bf 8} (1982-83), 16--20.
\bibitem{W} W. Wilczy\'nski, {\it Density topologies}, in: Handbook of Measure Theory (Edited by E. Pap), Elsevier, Amsterdam, 2002, 675--702.
\bibitem{WW} W. Wilczy\'nski, {\it A category analogue of the density topology non-homeomorphic with the $\mc I$-density topology}, Positivity {\bf 23} (2019), 469--484.
\end{thebibliography}
\end{document}